\documentclass{article}
\usepackage{amsthm}
\usepackage{fleqn}
\usepackage{array}
\usepackage[latin1]{inputenc}
\usepackage{amssymb}
\usepackage{amsfonts}
\usepackage{amsmath}

\begin{document}

\newtheorem{thm}{Theorem}[section]
\newtheorem{prop}{Proposition}[section]
\newtheorem{lem}{Lemma}[section]
\newtheorem{cor}{Corollary}[section]

\numberwithin{equation}{section}

\begin{center}
{\Large Some minimal hypersurfaces in the $(n+1)$-space

\vspace{1mm}

with $2m$-norm}

\vspace{5mm}

Makoto SAKAKI and Ryota TANAKA
\end{center}

\vspace{2mm}

{\bf Abstract.} We study translation minimal hypersurfaces and separable minimal hypersurfaces in the $(n+1)$-space with $2m$-norm. 

\vspace{2mm}

{\bf Mathematics Subject Classification.} 53A35, 53C42, 52A20, 52A21, 46B20

\vspace{2mm}

{\bf Keywords.} normed space, minimal hypersurface, translation hypersurface, separable hypersurface

\section{Introduction}

It is interesting to generalize differential geometry in Euclidean spaces to that in normed spaces, which was classically discussed in \cite{Bu} and \cite{G}.  Curves in normed spaces (generally, in gauge spaces) are studied in  \cite{BMSh}, \cite{BMS1} and \cite{BMS2}. Surfaces in normed $3$-spaces are studied in \cite{BMT1}, \cite{BMT2}, \cite{BMT3}, \cite{BMT4}. Rotational surfaces in a normed $3$-space with rotationally symmetric norm are discussed in \cite{S} and \cite{SY}. In the previous paper \cite{ST}, we consider some classes of minimal surfaces in the $3$-space with $2m$-norm. 

In this paper, generalizing \cite{ST}, we study translation minimal hypersurfaces and separable minimal hypersurfaces in the $(n+1)$-space with $2m$-norm. See \cite{Se} and \cite{CWW} for translation hypersurfaces and separable hypersurfaces in the Euclidean space. 

This paper is organized as follows. In Section 2, following and generalizing \cite{BMT2}, we give preliminaries on hypersurfaces in normed $(n+1)$-spaces. In Section 3, we compute the mean curvature of translation hypersurfaces in the $(n+1)$-space $({\mathbb R}^{n+1}, \|\cdot\|_{2m})$ with $2m$-norm.  In Section 4, translation minimal hypersurfaces in $({\mathbb R}^{n+1}, \|\cdot\|_{2m})$ are studied. In Section 5, we compute the mean curvature of separable hypersurfaces in $({\mathbb R}^{n+1}, \|\cdot\|_{2m})$. In Section 6, we discuss separable minimal hypersurfaces in $({\mathbb R}^{n+1}, \|\cdot\|_{2m})$ together with some examples.

\section{Preliminaries}

Let $({\mathbb R}^{n+1}, \|\cdot\|)$ be an $(n+1)$-dimensional normed space. The unit ball $B$ and the unit sphere $S$ are defined by
\[B = \{x \in {\mathbb R}^{n+1}; \| x\| \leq 1\}, \ \ \ \ S = \{x \in {\mathbb R}^{n+1}; \| x\| = 1\}. \]
We assume that $S$ is smooth and strictly convex, that is, $S$ is a smooth hypersurface in ${\mathbb R}^{n+1}$ and $S$ contains no line segment. 

Let $v$ be a non-zero vector in ${\mathbb R}^{n+1}$ and $\Pi$ be a hyperplane in ${\mathbb R}^{n+1}$. We say that $v$ is Birkhoff orthogonal to $\Pi$, denoted by $v \dashv_{B} \Pi$, if the tangent space of $S$ at $v/\|v\|$ is parallel to $\Pi$. 

Let $M$ be a hypersurface immersed in $({\mathbb R}^{n+1}, \|\cdot\|)$. Let $T_{p}M$ be the tangent space of $M$ at $p \in M$. Then there exists a vector $\eta(p) \in S$ such that $\eta(p) \dashv_{B} T_{p}M$. This gives a local smooth map $\eta: U \subset M \rightarrow S$, that is called the Birkhoff-Gauss map. It can be global if and only if $M$ is orientable. 

The mean curvature $H$  of $M$ is defined by
\[H = \frac{1}{n} \mbox{trace}(d\eta). \]
We say that $M$ is minimal if $H = 0$ identically.

\section{Translation hypersurfaces}

Let $({\mathbb R}^{n+1}, \|\cdot\|_{2m})$ be the $(n+1)$-space with $2m$-norm
\[\|x\|_{2m} = \left( \sum_{i=1}^{n+1} x_{i}^{2m} \right)^{\frac{1}{2m}}, \ \ \ x = (x_{1}, x_{2}, \cdots, x_{n+1}) \in {\mathbb R}^{n+1}, \]
where $m$ is a positive integer. Set
\[\Phi(x) := \sum_{i=1}^{n+1} x_{i}^{2m}. \]
The unit sphere $S$ is given by
\[S = \{x \in {\mathbb R}^{n+1} ; \Phi(x) = 1\}. \]

Let $M$ be a graph hypersurface in $({\mathbb R}^{n+1}, \|\cdot\|_{2m})$ given by
\[X(u_1, u_2, \cdots, u_n) = (u_1, u_2, \cdots, u_n, f(u_1, u_2, \cdots, u_n)) \]
for a smooth function $f(u_1, u_2, \cdots, u_n)$ with $f_{u_i} \neq 0$ for all $1 \leq i \leq n$. Then
\[X_{u_1} = (1, 0, \cdots, 0, f_{u_1}), \ X_{u_2} = (0, 1, 0, \cdots, 0, f_{u_2}), \ \cdots, \ X_{u_n} = (0, \cdots, 0, 1, f_{u_n}). \]
Set
\[\nu := (-f_{u_1}, -f_{u_2}, \cdots,  -f_{u_n}, 1) \]
which is orthogonal to $X_{u_i}$ ($1 \leq i \leq n$) with respect to the standard Euclidean metric. 

The Birkhoff-Gauss map $\eta = \eta(u_1, u_2, \cdots, u_n)$ is characterized by the condition
\[(\mbox{grad}(\Phi))_{\eta} = \left( \frac{\partial\Phi}{\partial x_1}(\eta), \frac{\partial\Phi}{\partial x_2}(\eta), \cdots, \frac{\partial\Phi}{\partial x_{n+1}}(\eta) \right) = \varphi\nu, \]
where $\varphi$ is a positive function. Then we can get
\[\eta = A^{-\frac{1}{2m}} \left( -(f_{u_1})^{\frac{1}{2m-1}}, -(f_{u_2})^{\frac{1}{2m-1}}, \cdots, -(f_{u_n})^{\frac{1}{2m-1}}, 1 \right) \]
where
\[A := 1+\sum_{i=1}^{n} (f_{u_i})^{\frac{2m}{2m-1}}. \]

In the following, we consider the case where $M$ is a translation hypersurface such that
\[f(u_1, u_2, \cdots, u_n) = f_{1}(u_1)+f_{2}(u_2)+\cdots +f_{n}(u_n), \]
where $f_{i}' \neq 0$ for all $1 \leq i \leq n$. Then
\[\eta = A^{-\frac{1}{2m}} \left( -(f_{1}')^{\frac{1}{2m-1}}, -(f_{2}')^{\frac{1}{2m-1}}, \cdots, -(f_{n}')^{\frac{1}{2m-1}}, 1 \right) \]
where
\[A = 1+\sum_{i=1}^{n} (f_{i}')^{\frac{2m}{2m-1}}. \]

We can compute that 
\[\eta_{u_j} = \sum_{k=1}^{n} \eta_{j}^{k} X_{u_k}, \ \ \ 1 \leq j \leq n \]
where
\[\eta_{j}^{j} = -\frac{A^{-\frac{2m+1}{2m}}}{2m-1} (f_{j}')^{-\frac{2m-2}{2m-1}} f_{j}'' \left\{ 1+\sum_{i\neq j} (f_{i}')^{\frac{2m}{2m-1}} \right\} \]
for $1 \leq j \leq n$, and
\[\eta_{j}^{k} = \frac{A^{-\frac{2m+1}{2m}}}{2m-1} (f_{j}')^{\frac{1}{2m-1}} f_{j}'' (f_{k}')^{\frac{1}{2m-1}} \]
for $j \neq k$ and $1 \leq j, k \leq n$. 

So the mean curvature $H$ is given by
\begin{eqnarray}
H = -\frac{A^{-\frac{2m+1}{2m}}}{n(2m-1)} \sum_{j=1}^{n} \left[ (f_{j}')^{-\frac{2m-2}{2m-1}} f_{j}'' \left\{ 1+\sum_{i\neq j} (f_{i}')^{\frac{2m}{2m-1}} \right\} \right]. 
\end{eqnarray}

\section{Translation minimal hypersurfaces}

Let $M$ be a translation hypersurface in $({\mathbb R}^{n+1}, \|\cdot\|_{2m})$ as in Section 3. 

\vspace{2mm}

(i) When $n = 3$, by (3.1), $M$ is minimal if and only if
\begin{eqnarray}
(f_{1}')^{-\frac{2m-2}{2m-1}} f_{1}'' \left\{ 1+(f_{2}')^{\frac{2m}{2m-1}}+(f_{3}')^{\frac{2m}{2m-1}} \right\} \nonumber 
\end{eqnarray}
\begin{eqnarray}
\hspace{5mm}+(f_{2}')^{-\frac{2m-2}{2m-1}} f_{2}'' \left\{ 1+(f_{1}')^{\frac{2m}{2m-1}}+(f_{3}')^{\frac{2m}{2m-1}} \right\} \nonumber 
\end{eqnarray}
\begin{eqnarray}
\hspace{5mm} +(f_{3}')^{-\frac{2m-2}{2m-1}} f_{3}'' \left\{ 1+(f_{1}')^{\frac{2m}{2m-1}}+(f_{2}')^{\frac{2m}{2m-1}} \right\} = 0. 
\end{eqnarray}

\vspace{1mm}

(i-1) We consider the case where $f_{1}'' f_{2}'' f_{3}'' \neq 0$. Differentiating (4.1) by $u_1$, we have
\[\left\{ -\frac{2m-2}{2m-1}(f_{1}')^{-\frac{4m-3}{2m-1}} (f_{1}'')^{2}+(f_{1}')^{-\frac{2m-2}{2m-1}}f_{1}^{(3)} \right\} \left\{ 1+(f_{2}')^{\frac{2m}{2m-1}}+(f_{3}')^{\frac{2m}{2m-1}} \right\} \]
\[
+\frac{2m}{2m-1}(f_{1}')^{\frac{1}{2m-1}}f_{1}'' \left\{ (f_{2}')^{-\frac{2m-2}{2m-1}} f_{2}'' +(f_{3}')^{-\frac{2m-2}{2m-1}} f_{3}'' \right\} = 0, \]
and
\begin{eqnarray}
(2m-2)\frac{f_{1}''}{(f_{1}')^{2}}-(2m-1)\frac{f_{1}^{(3)}}{f_{1}' f_{1}''} \nonumber
\end{eqnarray}
\begin{eqnarray}
\hspace{5mm} = 2m\cdot\frac{(f_{2}')^{-\frac{2m-2}{2m-1}} f_{2}'' +(f_{3}')^{-\frac{2m-2}{2m-1}} f_{3}''} {1+(f_{2}')^{\frac{2m}{2m-1}}+(f_{3}')^{\frac{2m}{2m-1}}} = c_1 
\end{eqnarray}
for a constant $c_1$. From the latter part of (4.2), we have
\begin{eqnarray}
2m(f_{2}')^{-\frac{2m-2}{2m-1}} f_{2}''-c_{1}(f_{2}')^{\frac{2m}{2m-1}} \nonumber
\end{eqnarray}
\begin{eqnarray}
\hspace{5mm} = c_{1}-\left\{ 2m(f_{3}')^{-\frac{2m-2}{2m-1}} f_{3}''-c_{1}(f_{3}')^{\frac{2m}{2m-1}} \right\} = c_{2} 
\end{eqnarray}
for a constant $c_2$. 

Analogously, differentiating (4.1) by $u_2$, we have
\begin{eqnarray}
(2m-2)\frac{f_{2}''}{(f_{2}')^{2}}-(2m-1)\frac{f_{2}^{(3)}}{f_{2}' f_{2}''} \nonumber
\end{eqnarray}
\begin{eqnarray}
\hspace{5mm} = 2m\cdot\frac{(f_{1}')^{-\frac{2m-2}{2m-1}} f_{1}'' +(f_{3}')^{-\frac{2m-2}{2m-1}} f_{3}''} {1+(f_{1}')^{\frac{2m}{2m-1}}+(f_{3}')^{\frac{2m}{2m-1}}} = c_3, 
\end{eqnarray}
and
\begin{eqnarray}
2m(f_{3}')^{-\frac{2m-2}{2m-1}} f_{3}''-c_{3}(f_{3}')^{\frac{2m}{2m-1}} \nonumber
\end{eqnarray}
\begin{eqnarray}
\hspace{5mm} = c_{3}-\left\{ 2m(f_{1}')^{-\frac{2m-2}{2m-1}} f_{1}''-c_{3}(f_{1}')^{\frac{2m}{2m-1}} \right\} = c_{4} 
\end{eqnarray}
for constants $c_3$ and $c_4$. Differentiating (4.1) by $u_3$, we have
\begin{eqnarray}
(2m-2)\frac{f_{3}''}{(f_{3}')^{2}}-(2m-1)\frac{f_{3}^{(3)}}{f_{3}' f_{3}''} \nonumber
\end{eqnarray}
\begin{eqnarray}
\hspace{5mm} = 2m\cdot\frac{(f_{1}')^{-\frac{2m-2}{2m-1}} f_{1}'' +(f_{2}')^{-\frac{2m-2}{2m-1}} f_{2}''} {1+(f_{1}')^{\frac{2m}{2m-1}}+(f_{2}')^{\frac{2m}{2m-1}}} = c_5, 
\end{eqnarray}
and
\begin{eqnarray}
2m(f_{1}')^{-\frac{2m-2}{2m-1}} f_{1}''-c_{5}(f_{1}')^{\frac{2m}{2m-1}} \nonumber
\end{eqnarray}
\begin{eqnarray}
\hspace{5mm} = c_{5}-\left\{ 2m(f_{2}')^{-\frac{2m-2}{2m-1}} f_{2}''-c_{5}(f_{2}')^{\frac{2m}{2m-1}} \right\} = c_{6} 
\end{eqnarray}
for constants $c_5$ and $c_6$. 

By (4.7) and (4.5), we have
\[2m(f_{1}')^{-\frac{2m-2}{2m-1}} f_{1}''-c_{5}(f_{1}')^{\frac{2m}{2m-1}} = c_{6}, \ \ \ 2m(f_{1}')^{-\frac{2m-2}{2m-1}} f_{1}''-c_{3}(f_{1}')^{\frac{2m}{2m-1}} = c_{3}-c_{4}. \]
If $c_3 \neq c_5$, then $f_{1}'$ is constant, which contradicts that $f_{1}'' \neq 0$. So $c_3 = c_5$ and $c_{3} = c_{4}+c_{6}$. Analogously, $c_1 = c_5$, $c_{5} = c_{2}+c_{6}$, $c_1 = c_3$ and  $c_{1} = c_{2}+c_{4}$. Thus we get
\[c_2 = c_4 = c_6 =: c_0, \ \ \ \ c_1 =c_3 = c_5 = 2c_0. \]
Hence, 
\begin{eqnarray}
2m(f_{i}')^{-\frac{2m-2}{2m-1}} f_{i}''-2c_{0}(f_{i}')^{\frac{2m}{2m-1}} = c_{0}, 
\end{eqnarray}
and
\begin{eqnarray}
f_{i}'' = \frac{c_0}{2m}\left\{ (f_{i}')^{\frac{2m-2}{2m-1}}+2(f_{i}')^{2} \right\} 
\end{eqnarray}
for $1 \leq i \leq 3$. 

On the other hand, from the former parts of (4.2), (4.4) and (4.6), we have
\begin{eqnarray}
(2m-2)\frac{f_{i}''}{(f_{i}')^{2}}-(2m-1)\frac{f_{i}^{(3)}}{f_{i}' f_{i}''} = 2c_{0}, \ \ \ 1 \leq i \leq 3, 
\end{eqnarray}
By (4.9) and (4.10) we find that $c_0 = 0$. Then $f_{i}'' = 0$ by (4.8), which is a contradiction. 

\vspace{1mm}

(i-2) If $f_{3}'' = 0$ for example, then $f_{3}(u_3) = au_{3}+b$ for constants $a$ and $b$. The equation (4.1) reduces to the $2$-dimensional case by a homothetical change of $f_i$, $1 \leq i \leq 2$, and $M$ is a cylinder over it.

\vspace{2mm}

(ii) When $n = 4$, by (3.1), $M$ is minimal if and only if
\begin{eqnarray}
(f_{1}')^{-\frac{2m-2}{2m-1}} f_{1}'' \left\{ 1+(f_{2}')^{\frac{2m}{2m-1}}+(f_{3}')^{\frac{2m}{2m-1}}+(f_{4}')^{\frac{2m}{2m-1}} \right\} \nonumber 
\end{eqnarray}
\begin{eqnarray}
\hspace{5mm} +(f_{2}')^{-\frac{2m-2}{2m-1}} f_{2}'' \left\{ 1+(f_{1}')^{\frac{2m}{2m-1}}+(f_{3}')^{\frac{2m}{2m-1}}+(f_{4}')^{\frac{2m}{2m-1}} \right\} \nonumber 
\end{eqnarray}
\begin{eqnarray}
\hspace{5mm} +(f_{3}')^{-\frac{2m-2}{2m-1}} f_{3}'' \left\{ 1+(f_{1}')^{\frac{2m}{2m-1}}+(f_{2}')^{\frac{2m}{2m-1}}+(f_{4}')^{\frac{2m}{2m-1}} \right\} \nonumber 
\end{eqnarray}
\begin{eqnarray}
\hspace{5mm}  +(f_{4}')^{-\frac{2m-2}{2m-1}} f_{4}'' \left\{ 1+(f_{1}')^{\frac{2m}{2m-1}}+(f_{2}')^{\frac{2m}{2m-1}}+(f_{3}')^{\frac{2m}{2m-1}} \right\} = 0. 
\end{eqnarray}

\vspace{1mm}

(ii-1) We consider the case where $f_{1}'' f_{2}'' f_{3}'' f_{4}'' \neq 0$. Differentiating (4.11) by $u_1$, we have
{\small 
\[\left\{ -\frac{2m-2}{2m-1}(f_{1}')^{-\frac{4m-3}{2m-1}} (f_{1}'')^{2}+(f_{1}')^{-\frac{2m-2}{2m-1}}f_{1}^{(3)} \right\} \left\{ 1+(f_{2}')^{\frac{2m}{2m-1}}+(f_{3}')^{\frac{2m}{2m-1}}+(f_{4}')^{\frac{2m}{2m-1}} \right\} \]
\[
+\frac{2m}{2m-1}(f_{1}')^{\frac{1}{2m-1}}f_{1}'' \left\{ (f_{2}')^{-\frac{2m-2}{2m-1}} f_{2}'' +(f_{3}')^{-\frac{2m-2}{2m-1}} f_{3}'' +(f_{4}')^{-\frac{2m-2}{2m-1}} f_{4}'' \right\} = 0, \] }
and
\begin{eqnarray}
(2m-2)\frac{f_{1}''}{(f_{1}')^{2}}-(2m-1)\frac{f_{1}^{(3)}}{f_{1}' f_{1}''} \nonumber
\end{eqnarray}
\begin{eqnarray}
\hspace{5mm} = 2m\cdot\frac{(f_{2}')^{-\frac{2m-2}{2m-1}} f_{2}'' +(f_{3}')^{-\frac{2m-2}{2m-1}} f_{3}'' +(f_{4}')^{-\frac{2m-2}{2m-1}} f_{4}''} {1+(f_{2}')^{\frac{2m}{2m-1}}+(f_{3}')^{\frac{2m}{2m-1}}+(f_{4}')^{\frac{2m}{2m-1}}} = c_1 
\end{eqnarray}
for a constant $c_1$. From the latter part of (4.12), we have
\begin{eqnarray}
2m(f_{2}')^{-\frac{2m-2}{2m-1}} f_{2}''-c_{1}(f_{2}')^{\frac{2m}{2m-1}} = c_{1}-\left\{ 2m(f_{3}')^{-\frac{2m-2}{2m-1}} f_{3}''-c_{1}(f_{3}')^{\frac{2m}{2m-1}} \right\} \nonumber
\end{eqnarray}
\begin{eqnarray}
\hspace{1cm} -\left\{ 2m(f_{4}')^{-\frac{2m-2}{2m-1}} f_{4}''-c_{1}(f_{4}')^{\frac{2m}{2m-1}} \right\} = c_{2} 
\end{eqnarray}
for a constant $c_2$. Here we should notice that
\[2m(f_{3}')^{-\frac{2m-2}{2m-1}} f_{3}''-c_{1}(f_{3}')^{\frac{2m}{2m-1}}, \ \ \ \ 2m(f_{4}')^{-\frac{2m-2}{2m-1}} f_{4}''-c_{1}(f_{4}')^{\frac{2m}{2m-1}} \]
are also constant. 

Analogousely, we have
\begin{eqnarray}
(2m-2)\frac{f_{2}''}{(f_{2}')^{2}}-(2m-1)\frac{f_{2}^{(3)}}{f_{2}' f_{2}''} \nonumber
\end{eqnarray}
\begin{eqnarray}
\hspace{5mm} = 2m\cdot\frac{(f_{1}')^{-\frac{2m-2}{2m-1}} f_{1}'' +(f_{3}')^{-\frac{2m-2}{2m-1}} f_{3}'' +(f_{4}')^{-\frac{2m-2}{2m-1}} f_{4}''} {1+(f_{1}')^{\frac{2m}{2m-1}}+(f_{3}')^{\frac{2m}{2m-1}}+(f_{4}')^{\frac{2m}{2m-1}}} = c_3, 
\end{eqnarray}
\begin{eqnarray}
2m(f_{3}')^{-\frac{2m-2}{2m-1}} f_{3}''-c_{3}(f_{3}')^{\frac{2m}{2m-1}} = c_{3}-\left\{ 2m(f_{4}')^{-\frac{2m-2}{2m-1}} f_{4}''-c_{3}(f_{4}')^{\frac{2m}{2m-1}} \right\} \nonumber
\end{eqnarray}
\begin{eqnarray}
\hspace{1cm} -\left\{ 2m(f_{1}')^{-\frac{2m-2}{2m-1}} f_{1}''-c_{3}(f_{1}')^{\frac{2m}{2m-1}} \right\} = c_{4}, 
\end{eqnarray}
\begin{eqnarray}
(2m-2)\frac{f_{3}''}{(f_{3}')^{2}}-(2m-1)\frac{f_{3}^{(3)}}{f_{3}' f_{3}''} \nonumber
\end{eqnarray}
\begin{eqnarray}
\hspace{5mm} = 2m\cdot\frac{(f_{1}')^{-\frac{2m-2}{2m-1}} f_{1}'' +(f_{2}')^{-\frac{2m-2}{2m-1}} f_{2}'' +(f_{4}')^{-\frac{2m-2}{2m-1}} f_{4}''} {1+(f_{1}')^{\frac{2m}{2m-1}}+(f_{2}')^{\frac{2m}{2m-1}}+(f_{4}')^{\frac{2m}{2m-1}}} = c_5, 
\end{eqnarray}
\begin{eqnarray}
2m(f_{4}')^{-\frac{2m-2}{2m-1}} f_{4}''-c_{5}(f_{4}')^{\frac{2m}{2m-1}} = c_{5}-\left\{ 2m(f_{1}')^{-\frac{2m-2}{2m-1}} f_{1}''-c_{5}(f_{1}')^{\frac{2m}{2m-1}} \right\} \nonumber
\end{eqnarray}
\begin{eqnarray}
\hspace{1cm} -\left\{ 2m(f_{2}')^{-\frac{2m-2}{2m-1}} f_{2}''-c_{5}(f_{2}')^{\frac{2m}{2m-1}} \right\} = c_{6}, 
\end{eqnarray}
\begin{eqnarray}
(2m-2)\frac{f_{4}''}{(f_{4}')^{2}}-(2m-1)\frac{f_{4}^{(3)}}{f_{4}' f_{4}''} \nonumber
\end{eqnarray}
\begin{eqnarray}
\hspace{5mm} = 2m\cdot\frac{(f_{1}')^{-\frac{2m-2}{2m-1}} f_{1}'' +(f_{2}')^{-\frac{2m-2}{2m-1}} f_{2}'' +(f_{3}')^{-\frac{2m-2}{2m-1}} f_{3}''} {1+(f_{1}')^{\frac{2m}{2m-1}}+(f_{2}')^{\frac{2m}{2m-1}}+(f_{3}')^{\frac{2m}{2m-1}}} = c_7, 
\end{eqnarray}
\begin{eqnarray}
2m(f_{1}')^{-\frac{2m-2}{2m-1}} f_{1}''-c_{7}(f_{1}')^{\frac{2m}{2m-1}} = c_{7}-\left\{ 2m(f_{2}')^{-\frac{2m-2}{2m-1}} f_{2}''-c_{7}(f_{2}')^{\frac{2m}{2m-1}} \right\} \nonumber
\end{eqnarray}
\begin{eqnarray}
\hspace{1cm} -\left\{ 2m(f_{3}')^{-\frac{2m-2}{2m-1}} f_{3}''-c_{7}(f_{3}')^{\frac{2m}{2m-1}} \right\} = c_{8} 
\end{eqnarray}
for constants $c_i$, $3 \leq i \leq 8$. 

By (4.13), (4.15), (4.17) and (4.19), as in (i-1), we can find that
\[c_2 = c_4 = c_6 = c_8 =: c_0, \ \ \ \ c_1 =c_3 = c_5 = c_7 = 3c_0. \]
So
\begin{eqnarray}
2m(f_{i}')^{-\frac{2m-2}{2m-1}} f_{i}''-3c_{0}(f_{i}')^{\frac{2m}{2m-1}} = c_{0}, 
\end{eqnarray}
and
\begin{eqnarray}
f_{i}'' = \frac{c_0}{2m}\left\{ (f_{i}')^{\frac{2m-2}{2m-1}}+3(f_{i}')^{2} \right\} 
\end{eqnarray}
for $1 \leq i \leq 4$. 

From the former parts of (4.12), (4.14), (4.16) and (4.18), we have
\begin{eqnarray}
(2m-2)\frac{f_{i}''}{(f_{i}')^{2}}-(2m-1)\frac{f_{i}^{(3)}}{f_{i}' f_{i}''} = 3c_{0}, \ \ \ 1 \leq i \leq 4. 
\end{eqnarray}
By (4.21) and (4.22) we find that $c_0 = 0$. Then $f_{i}'' = 0$ by (4.20), which is a contradiction. 

\vspace{1mm}

(ii-2) If $f_{4}'' = 0$ for example, then $f_{4}(u_4) = au_{4}+b$ for constants $a, b$, and the equation (4.11) is equivalent to (4.1) by a homothetical change of $f_i$, $1 \leq i \leq 3$. Then by (i),  $f_i'' = 0$ for some $1 \leq i \leq 3$. Thus the discussion reduces to the $2$-dimensional case, and $M$ is a cylinder over it.

\vspace{2mm}

(iii) In general, by (3.1), $M$ is minimal if and only if
\begin{eqnarray}
\hspace{5mm} \sum_{j=1}^{n} \left[ (f_{j}')^{-\frac{2m-2}{2m-1}} f_{j}'' \left\{ 1+\sum_{i\neq j} (f_{i}')^{\frac{2m}{2m-1}} \right\} \right] = 0. 
\end{eqnarray}
Inductively, we can see that at least $n-2$ of $f_{i}$, $1 \leq i \leq n$ are linear functions, and the equation (4.23) reduces to the $2$-dimensional case. Thus we obtain the following: 

\begin{thm}
Let $M$ be a translation hypersurface in $({\mathbb R}^{n+1}, \|\cdot\|_{2m})$ as above. Then $M$ is minimal if and only if it is either a hyperplane or a cylinder over a non-planar translation minimal surface in a $3$-space. 
\end{thm}

{\bf Remark.} The above mentioned non-planar translation minimal surface in a $3$-space is given in \cite{ST}. Here some homothetical change is considered.

\section{Separable hypersurfaces}

Let $M$ be an implicit hypersurface in $({\mathbb R}^{n+1}, \|\cdot\|_{2m})$ which is given by
\[f_{1}(x_{1})+f_{2}(x_{2})+\cdots +f_{n+1}(x_{n+1}) = 0 \]
for smooth functions $f_{i}(x_{i})$, $1 \leq i \leq n+1$ with 
\[\left( f_{1}'(x_{1}), f_{2}'(x_{2}), \cdots, f_{n+1}'(x_{n+1}) \right) \neq (0, 0, \cdots, 0) \]
on $M$. Such a hypersurface in called a separable hypersurface. Translation hypersurfaces are examples of separable hypersurfaces. 

We assume that $M$ is not a translation hypersurface. Then we may assume that 
\[\prod_{i=1}^{n+1} f_{i}'(x_{i}) \neq 0, \ \ \ \ \prod_{i=1}^{n+1} f_{i}''(x_{i}) \neq 0. \]

The Birkhoff-Gauss map $\eta$ is given by the condition
\[(\mbox{grad}(\Phi))_{\eta} = \varphi (f_{1}', f_{2}', \cdots, f_{n+1}'), \]
where $\varphi$ is a positive function. Then we have
\[\eta = A^{-\frac{1}{2m}} \left( (f_{1}')^{\frac{1}{2m-1}}, (f_{2}')^{\frac{1}{2m-1}}, \cdots, (f_{n+1}')^{\frac{1}{2m-1}} \right) \]
where
\[A := \sum_{i=1}^{n+1} (f_{i}')^{\frac{2m}{2m-1}}. \]

Since $f_{n+1}'(x_{n+1}) \neq 0$, $M$ can be locally expressed as $x_{n+1} = x_{n+1}(x_1, x_2, \cdots, x_n)$ and
\[\frac{\partial x_{n+1}}{\partial x_j} = -\frac{f_{j}'}{f_{n+1}'}, \ \ \ 1 \leq j \leq n. \]
So $M$ is locally parametrized as 
\[p(x_1, x_2, \cdots, x_n) = (x_1, x_2, \cdots, x_n, x_{n+1}(x_1, x_2, \cdots, x_n)) \]
and
\[p_{x_j} = \left( 0, \cdots, 0, 1, 0, \cdots, 0, -\frac{f_{j}'}{f_{n+1}'} \right), \]
where the $j$-th component is $1$. 

With respect to the local coordinates $x_1, x_2, \cdots, x_n$, we can compute that 
\[\eta_{x_j} = \sum_{k=1}^{n}\eta_{j}^{k} p_{x_k}, \]
where
\[\eta_{j}^{j} = \frac{A^{-\frac{2m+1}{2m}}}{2m-1} \left\{ (f_j')^{\frac{2m}{2m-1}} (f_{n+1}')^{-\frac{2m-2}{2m-1}} f_{n+1}''+ (f_{j}')^{-\frac{2m-2}{2m-1}} f_{j}'' \sum_{i\neq j}^{n+1} (f_{i}')^{\frac{2m}{2m-1}} \right\} \]
for $1 \leq j \leq n$, and
\[\eta_{j}^{k} = \frac{A^{-\frac{2m+1}{2m}}}{2m-1} (f_{k}')^{\frac{1}{2m-1}} \left\{ f_{j}' (f_{n+1}')^{-\frac{2m-2}{2m-1}} f_{n+1}''-(f_{j}')^{\frac{1}{2m-1}} f_{j}'' \right\} \]
for $j \neq k$ and $1 \leq j, k \leq n$. 

So the mean curvature $H$ is given by
\begin{eqnarray}
H = \frac{A^{-\frac{2m+1}{2m}}}{n(2m-1)} \sum_{j=1}^{n+1} \left\{ (f_{j}')^{-\frac{2m-2}{2m-1}} f_{j}'' \sum_{i\neq j}^{n+1} (f_{i}')^{\frac{2m}{2m-1}} \right\}. 
\end{eqnarray}

\section{Separable minimal hypersurfaces}

Let $M$ be a separable hypersurface in $({\mathbb R}^{n+1}, \|\cdot\|_{2m})$ as in Section 5. By (5.1), $M$ is minimal if and only if
\begin{eqnarray}
\sum_{j=1}^{n+1} \left\{ (f_{j}')^{-\frac{2m-2}{2m-1}} f_{j}'' \sum_{i\neq j}^{n+1} (f_{i}')^{\frac{2m}{2m-1}} \right\} \nonumber 
\end{eqnarray}
\begin{eqnarray}
\hspace{5mm} = \sum_{j=1}^{n+1} \left\{ (f_{j}')^{-\frac{2m-2}{2m-1}} f_{j}'' \left( A-(f_{j}')^{\frac{2m}{2m-1}} \right) \right\} = 0. 
\end{eqnarray}

We continue the discussion as in \cite{KL}. We introduce new variables $u_1, u_2,..., u_n, u_{n+1}$ as
\[u_i = f_{i}(x_i), \ \ \ 1 \leq i \leq n+1. \]
Set
\[X_{i}(u_i) = (f_{i}')^{\frac{2m}{2m-1}}, \ \ \ 1 \leq i \leq n+1. \]
Then $A = \sum_{i=1}^{n+1} X_{i}$ and 
\[X_{i}' = \frac{dX_i}{du_i} = \frac{dX_i}{dx_i} \frac{dx_i}{du_i} = \frac{2m}{2m-1} (f_{i}')^{-\frac{2m-2}{2m-1}} f_{i}'' \neq 0, \ \ \ 1 \leq i \leq n+1. \]
We rewrite the equation (6.1) as
\begin{eqnarray}
\hspace{5mm} \sum_{j=1}^{n+1} \left( X_{j}' \sum_{i\neq j}^{n+1}X_{i} \right) =  \sum_{j=1}^{n+1} X_{j}' (A-X_{j}) = 0 
\end{eqnarray}
for all values $u_1, u_2, \cdots, u_{n+1}$ with $u_{1}+u_{2}+\cdots +u_{n+1} = 0$. 

Since
\[f_{i}' = \frac{df_i}{dx_i} = \frac{du_i}{dx_i} = \pm(X_{i}(u_i))^{\frac{2m-1}{2m}}, \]
the hypersurface can be represented as follows: 

\begin{thm}
Under the noation as above, the separable minimal hypersurface $M$ can be represented as 
\[x_{i} = \pm\int^{u_i} (X_{i}(u_i))^{-\frac{2m-1}{2m}} du_{i}, \ \ \ 1 \leq i \leq n, \]
\[x_{n+1} = \pm\int^{u_{n+1}=-u_{1}-u_{2}-\cdots -u_{n}} (X_{n+1}(u_{n+1}))^{-\frac{2m-1}{2m}} du_{n+1}, \]
where $X_i$ ($1 \leq i \leq n+1)$ should be positive and satisfy (6.2). 
\end{thm}

(i) We consider the case where $n = 3$ and $X_{i}(u_i) = p_{i}+q_{i}u_{i}$ for $1 \leq i \leq 4$, with $q_i \neq 0$ for all $i$. Then the equation (6.2) becomes the following identity: 
\[P_{0}+P_{1}u_{1}+P_{2}u_{2}+P_{3}u_{3} = 0. \]
The coefficients should be zero, so that
{\small 
\begin{eqnarray}
q_{1}(p_{2}+p_{3}+p_{4})+q_{2}(p_{1}+p_{3}+p_{4})+q_{3}(p_{1}+p_{2}+p_{4})+q_{4}(p_{1}+p_{2}+p_{3}) = 0, \nonumber 
\end{eqnarray} 
\begin{eqnarray}
(q_{2}+q_{3}) (q_{1}-q_{4}) = 0, \ \ \ (q_{1}+q_{3}) (q_{2}-q_{4}) = 0, \ \ \ (q_{1}+q_{2}) (q_{3}-q_{4}) = 0. 
\end{eqnarray} }
The domain is given by
\[p_{1}+q_{1}u_{1} > 0, \ \ \ p_{2}+q_{2}u_{2} > 0, \ \ \ p_{3}+q_{3}u_{3} > 0, \]
\[p_{4}-q_{4}(u_{1}+u_{2}+u_{3}) > 0. \]
By (6.3), it suffices to consider the following two cases: 

\ (i-1) $q_1 = q_2 = q_3 = q_4$ and $p_{1}+p_{2}+p_{3}+p_{4} = 0$; or

\ (i-2) $q_1 = -q_2 = -q_3 = q_4$ and $-p_{1}+p_{2}+p_{3}-p_{4} = 0$. 

\vspace{2mm}

(i-1) In this case, the domain is given by
\[p_{1}+q_{1}u_{1} > 0, \ \ \ p_{2}+q_{1}u_{2} > 0, \ \ \ p_{3}+q_{1}u_{3} > 0, \ \ \ p_{4}-q_{1}(u_{1}+u_{2}+u_{3}) > 0. \]
Adding them we have a contradiction. So this case (i-1) does not occur. 

\vspace{1mm}

(i-2) In this case, we have
\[x_{1} = \pm\int^{u_1} (p_{1}+q_{1}u_{1})^{-\frac{2m-1}{2m}} du_{1} = \pm\frac{2m}{q_1} (p_{1}+q_{1}u_{1})^{\frac{1}{2m}}, \]
\[x_{2} = \pm\int^{u_2} (p_{2}-q_{1}u_{2})^{-\frac{2m-1}{2m}} du_{2} = \mp\frac{2m}{q_1} (p_{2}-q_{1}u_{2})^{\frac{1}{2m}}, \]
\[x_{3} = \pm\int^{u_3} (p_{3}-q_{1}u_{3})^{-\frac{2m-1}{2m}} du_{3} = \mp\frac{2m}{q_1} (p_{3}-q_{1}u_{3})^{\frac{1}{2m}}, \]
\[x_{4} =  \pm\int^{u_{4}=-u_{1}-u_{2}-u_{3}}  (p_{4}+q_{1}u_{4})^{-\frac{2m-1}{2m}} du_{4} = \pm\frac{2m}{q_1} \left\{ p_{4}-q_{1}(u_{1}+u_{2}+u_{3}) \right\}^{\frac{1}{2m}}. \]
We rewrite them as
\[u_{1} = \frac{1}{q_1} \left\{ \left( \frac{q_{1}x_{1}}{2m} \right)^{2m}-p_{1} \right\} = f_{1}(x_{1}), \]
\[u_{2} = \frac{1}{q_1} \left\{ p_{2}-\left( \frac{q_{1}x_{2}}{2m} \right)^{2m} \right\} = f_{2}(x_{2}), \]
\[u_{3} = \frac{1}{q_1} \left\{ p_{3}-\left( \frac{q_{1}x_{3}}{2m} \right)^{2m} \right\} = f_{3}(x_{3}), \]
\[u_{4} = -u_{1}-u_{2}-u_{3} = \frac{1}{q_1} \left\{ \left( \frac{q_{1}x_{4}}{2m} \right)^{2m}-p_{4} \right\} = f_{4}(x_{4}). \]
Since $-p_{1}+p_{2}+p_{3}-p_{4} = 0$, as a separable hypersurface, we get
\[x_{1}^{2m}-x_{2}^{2m}-x_{3}^{2m}+x_{4}^{2m} = 0. \]

\vspace{2mm}

{\bf Example 6.1.} Set
\[p_1 = p_2 = p_3 = p_4 = 1, \ \ \ q_1 = -q_2 = -q_3 = q_4 = 1, \]
which is in the case (i-2). The domain is given by
\[1+u_{1} > 0, \ \ \ 1-u_{2} > 0, \ \ \ 1-u_{3} > 0, \ \ \ 1-u_{1}-u_{2}-u_{3} > 0, \]
which is non-empty. 

\vspace{2mm}

As a generalization, when $n = 2r-1$, we have the following: 

\vspace{2mm}

{\bf Example 6.2.} A separable hypersurface in $({\mathbb R}^{2r}, \|\cdot\|_{2m})$, $r \geq 2$, given by
\[\sum_{i=1}^{r}x_{i}^{2m}-\sum_{i=1}^{r}x_{r+i}^{2m} = 0 \]
is minimal, which can be shown by using (6.1).

\vspace{2mm}

(ii) We consider the case where $n = 3$ and $X_{i}(u_i) = p_{i}+q_{i}u_{i}+r_{i}u_{i}^{2}$ for $1 \leq i \leq 4$. Here $r_i \neq 0$ for some $i$, and $q_j \neq 0$ if $r_j = 0$. Then the equation (6.2) becomes the following identity: 
\[P_{0}+P_{1}u_{1}+P_{2}u_{2}+P_{3}u_{3}+P_{11}u_{1}^{2}+P_{22}u_{2}^{2}+P_{33}u_{3}^{2}+P_{12}u_{1}u_{2}+P_{13}u_{1}u_{3}+P_{23}u_{2}u_{3} \]
\[+P_{112}(u_{1}^{2}u_{2}+u_{1}u_{2}^{2})+P_{113}(u_{1}^{2}u_{3}+u_{1}u_{3}^{2})+P_{223}(u_{2}^{2}u_{3}+u_{2}u_{3}^{2})+P_{123}u_{1}u_{2}u_{3} = 0. \]
The coefficients are zero and we get
\[q_{1}(p_{2}+p_{3}+p_{4})+q_{2}(p_{1}+p_{3}+p_{4})+q_{3}(p_{1}+p_{2}+p_{4})+q_{4}(p_{1}+p_{2}+p_{3}) = 0, \]
\[(q_{2}+q_{3}) (q_{1}-q_{4})+2r_{1}(p_{2}+p_{3}+p_{4})-2r_{4}(p_{1}+p_{2}+p_{3}) = 0, \]
\[(q_{1}+q_{3}) (q_{2}-q_{4})+2r_{2}(p_{1}+p_{3}+p_{4})-2r_{4}(p_{1}+p_{2}+p_{3}) = 0, \]
\[(q_{1}+q_{2}) (q_{3}-q_{4})+2r_{3}(p_{1}+p_{2}+p_{4})-2r_{4}(p_{1}+p_{2}+p_{3}) = 0, \]
\[(q_{2}+q_{3})(r_{1}+r_{4})-q_{1}r_{4}-q_{4}r_{1} = 0, \]
\[(q_{1}+q_{3})(r_{2}+r_{4})-q_{2}r_{4}-q_{4}r_{2} = 0, \]
\[(q_{1}+q_{2})(r_{3}+r_{4})-q_{3}r_{4}-q_{4}r_{3} = 0, \]
\[(q_{2}-q_{4})r_{1}+(q_{1}-q_{4})r_{2}+q_{3}r_{4} = 0, \]
\[(q_{3}-q_{4})r_{1}+(q_{1}-q_{4})r_{3}+q_{2}r_{4} = 0, \]
\[(q_{3}-q_{4})r_{2}+(q_{2}-q_{4})r_{3}+q_{1}r_{4} = 0, \]
\[r_{1}r_{2}+r_{1}r_{4}+r_{2}r_{4} = 0, \ \ \ r_{1}r_{3}+r_{1}r_{4}+r_{3}r_{4} = 0, \]
\[r_{2}r_{3}+r_{2}r_{4}+r_{3}r_{4} = 0, \ \ \ r_{4}(r_{1}+r_{2}+r_{3}) = 0. \]

\vspace{1mm}

(ii-1) The case where $r_4 = 0$. We may assume that $r_1 \neq 0$. Then we can see that $r_2 = r_3 = 0$ and $q_2 = q_3 = q_4 = 0$, which is a contradiction. 

\vspace{1mm}

(ii-2) The case where $r_4 \neq 0$. Then
\[r_{1}+r_{2}+r_{3} = 0, \ \ \ r_{1}r_{2} = r_{3}r_{4}, \ \ \ r_{1}r_{3} = r_{2}r_{4}, \ \ \ r_{2}r_{3} = r_{1}r_{4}. \]
Since
\[r_{1}r_{2}r_{3} = r_{1}^{2}r_{4} = r_{2}^{2}r_{4} = r_{3}^{2}r_{4}, \]
we have $|r_1| = |r_2| = |r_3|$. Combining with $r_{1}+r_{2}+r_{3} = 0$, we see that $r_1 = r_2 = r_3 = 0$. Then we find that $q_1 = q_2 = q_3 = 0$, which is a contradiction. 

\vspace{1mm}

Therefore, the case (ii) does not occur.

\vspace{2mm}

(iii) We consider the case where $n = 4$ and $X_{i}(u_i) = p_{i}+q_{i}u_{i}$ for $1 \leq i \leq 5$, with $q_i \neq 0$ for all $i$. Then (6.2) becomes
\[P_{0}+P_{1}u_{1}+P_{2}u_{2}+P_{3}u_{3}+P_{4}u_{4} = 0, \]
so that
\begin{eqnarray}
q_{1}(p_{2}+p_{3}+p_{4}+p_{5})+q_{2}(p_{1}+p_{3}+p_{4}+p_{5})+q_{3}(p_{1}+p_{2}+p_{4}+p_{5}) \nonumber
\end{eqnarray}
\begin{eqnarray}
\hspace{5mm} +q_{4}(p_{1}+p_{2}+p_{3}+p_{5})+q_{5}(p_{1}+p_{2}+p_{3}+p_{4}) = 0,  \nonumber
\end{eqnarray}
\begin{eqnarray}
(q_{2}+q_{3}+q_{4})(q_{1}-q_{5}) = 0, \ \ \ \ (q_{1}+q_{3}+q_{4})(q_{2}-q_{5}) = 0,  \nonumber
\end{eqnarray}
\begin{eqnarray}
(q_{1}+q_{2}+q_{4})(q_{3}-q_{5}) = 0, \ \ \ \ (q_{1}+q_{2}+q_{3})(q_{4}-q_{5}) = 0. 
\end{eqnarray}
By (6.4), it suffices to consider the following cases: 

\ (iii-1) $q_1 = q_2 = q_3 = q_4 = q_5$; 

\ (iii-2) $q_1 = q_2 = q_5$ and $q_{1}+q_{2}+q_{3} = q_{1}+q_{2}+q_{4} = 0$; or

\ (iii-3) $q_1 = q_5$ and $q_{1}+q_{2}+q_{3} = q_{1}+q_{2}+q_{4} = q_{1}+q_{3}+q_{4} = 0$. 

\vspace{2mm}

(iii-1) In this case, we have $p_{1}+p_{2}+p_{3}+p_{4}+p_{5} = 0$. As in the case (i-1), we find that the domain is empty. 

\vspace{1mm}

(iii-2) In this case, we have
\[q_{3} = q_{4} = -2q_{1} = -2q_{2} = -2q_{5}, \ \ \ -2p_{1}-2p_{2}+p_{3}+p_{4}-2p_{5} = 0. \]
Then
\[x_{1} = \pm\int^{u_1} (p_{1}+q_{1}u_{1})^{-\frac{2m-1}{2m}} du_{1} = \pm\frac{2m}{q_1} (p_{1}+q_{1}u_{1})^{\frac{1}{2m}}, \]
\[x_{2} = \pm\int^{u_2} (p_{2}+q_{1}u_{2})^{-\frac{2m-1}{2m}} du_{2} = \pm\frac{2m}{q_1} (p_{2}+q_{1}u_{2})^{\frac{1}{2m}}, \]
\[x_{3} = \pm\int^{u_3} (p_{3}-2q_{1}u_{3})^{-\frac{2m-1}{2m}} du_{3} = \mp\frac{m}{q_1} (p_{3}-2q_{1}u_{3})^{\frac{1}{2m}}, \]
\[x_{4} = \pm\int^{u_4} (p_{4}-2q_{1}u_{4})^{-\frac{2m-1}{2m}} du_{4} = \mp\frac{m}{q_1} (p_{4}-2q_{1}u_{4})^{\frac{1}{2m}}, \]
\[x_{5} = \pm\int^{u_{5}=-u_{1}-u_{2}-u_{3}-u_{4}}  (p_{5}+q_{1}u_{5})^{-\frac{2m-1}{2m}} du_{5} \]
\[= \pm\frac{2m}{q_1} \left\{ p_{5}-q_{1}(u_{1}+u_{2}+u_{3}+u_{4}) \right\}^{\frac{1}{2m}}, \]
which are rewritten as
\[u_{1} = \frac{1}{q_1} \left\{ \left( \frac{q_{1}x_{1}}{2m} \right)^{2m}-p_{1} \right\} = f_{1}(x_{1}), \]
\[u_{2} = \frac{1}{q_1} \left\{ \left( \frac{q_{1}x_{2}}{2m} \right)^{2m}-p_{2} \right\} = f_{2}(x_{2}), \]
\[u_{3} = \frac{1}{2q_1} \left\{ p_{3}-\left( \frac{q_{1}x_{3}}{m} \right)^{2m} \right\} = f_{3}(x_{3}), \]
\[u_{4} = \frac{1}{2q_1} \left\{ p_{4}-\left( \frac{q_{1}x_{4}}{m} \right)^{2m} \right\} = f_{4}(x_{4}), \]
\[u_{5} = -u_{1}-u_{2}-u_{3}-u_{4} = \frac{1}{q_1} \left\{ \left( \frac{q_{1}x_{5}}{2m} \right)^{2m}-p_{5} \right\} = f_{5}(x_{5}). \]
Noticing that
\[-2p_{1}-2p_{2}+p_{3}+p_{4}-2p_{5} = 0, \]
as a separable hypersurface, we get
\[x_{1}^{2m}+x_{2}^{2m}-2^{2m-1}x_{3}^{2m}-2^{2m-1}x_{4}^{2m}+x_{5}^{2m} = 0. \]

\vspace{2mm}

{\bf Example 6.3.} Set
\[p_{1} = p_{2} = p_{5} = 1, \ \ \ p_{3} = p_{4} = 3, \ \ \ q_{1} = q_{2} = q_{5} = 1, \ \ \ q_{3} = q_{4} = -2, \]
which is in the case (iii-2). The domain is given by
\[1+u_1 > 0, \ \ \ 1+u_2 > 0, \ \ \ 3-2u_3 > 0, \ \ \ 3-2u_4 > 0, \]
\[1-u_{1}-u_{2}-u_{3}-u_{4} > 0, \]
which is non-empty. 

\vspace{2mm}

(iii-3) In this case, we have
\[q_1 = q_5 = -2q_2 = -2q_3 = -2q_4, \ \ \ \ p_{1}-2p_{2}-2p_{3}-2p_{4}+p_{5} = 0. \]
This case can be identified with (iii-2). 

\vspace{2mm}

As a generalizaation, when $n = 2r$, we have the following: 

\vspace{2mm}

{\bf Example 6.4.} A separable hypersurface in $({\mathbb R}^{2r+1}, \|\cdot\|_{2m})$, $r \geq 2$, given by
\[-\left( \frac{r}{r-1} \right)^{2m-1} \sum_{i=1}^{r} x_{i}^{2m}+\sum_{i=1}^{r+1} x_{r+i}^{2m} = 0 \]
is minimal, which can be shown by (6.1).

\vspace{2mm}

(iv) The case where $n = 3$ and $X_{i}(u_i) = q_{i}e^{u_i}+r_{i}e^{-u_i}$ for $1 \leq i \leq 4$. Then (6.2) becomes
\begin{eqnarray}
2(q_{1}q_{2}-r_{3}r_{4})e^{u_{1}+u_{2}}+2(q_{3}q_{4}-r_{1}r_{2})e^{-u_{1}-u_{2}} \nonumber 
\end{eqnarray}
\begin{eqnarray}
\hspace{3mm} +2(q_{1}q_{3}-r_{2}r_{4})e^{u_{1}+u_{3}}+2(q_{2}q_{4}-r_{1}r_{3})e^{-u_{1}-u_{3}} \nonumber 
\end{eqnarray}
\begin{eqnarray}
\hspace{3mm} +2(q_{2}q_{3}-r_{1}r_{4})e^{u_{2}+u_{3}}+2(q_{1}q_{4}-r_{2}r_{3})e^{-u_{2}-u_{3}} = 0. 
\end{eqnarray}

\vspace{2mm}

{\bf Example 6.5.} Set $q_i = 1$, $r_i = 1$ for $1 \leq i \leq 4$. Then (6.5) is valid and
\[X_1 = e^{u_1}+e^{-u_1}, \ \ X_2 = e^{u_2}+e^{-u_2}, \ \ X_3 = e^{u_3}+e^{-u_3}, \ \ X_4 = e^{u_4}+e^{-u_4}. \]
The domain is ${\mathbb R}^3 = \{(u_{1}, u_{2}, u_{3}); u_{i} \in {\mathbb R}, 1 \leq i \leq 3\}$, with $u_{4} = -u_{1}-u_{2}-u_{3}$. 

\vspace{2mm}

{\bf Example 6.6.} Set 
\[q_1= 1, \ r_1 = 0, \ q_2 = 0, \ r_2 = 1, \ q_3 = 0, \ r_3 = 1, \ q_4 = 1, \ r_4 = 0. \]
Then (6.5) holds and
\[X_1 = e^{u_1}, \ \ \ X_2 = e^{-u_2}, \ \ \ X_3 = e^{-u_3}, \ \ \ X_4 = e^{u_4}. \]
We have
\[x_1 = \pm\int^{u_1} e^{-\frac{2m-1}{2m}u_1} du_1 = \mp\frac{2m}{2m-1} e^{-\frac{2m-1}{2m}u_1}, \]
\[x_2 = \pm\int^{u_2} e^{\frac{2m-1}{2m}u_2} du_2 = \pm\frac{2m}{2m-1} e^{\frac{2m-1}{2m}u_2}, \]
\[x_3 = \pm\int^{u_3} e^{\frac{2m-1}{2m}u_3} du_3 = \pm\frac{2m}{2m-1} e^{\frac{2m-1}{2m}u_3}, \]
\[x_4 = \pm\int^{u_{4}=-u_{1}-u_{2}-u_{3}} e^{-\frac{2m-1}{2m}u_4} du_4 = \mp\frac{2m}{2m-1} e^{\frac{2m-1}{2m}(u_{1}+u_{2}+u_{3})}, \]
and rewrite them as
\[u_1 = -\frac{2m}{2m-1} \log{ \left( \frac{2m-1}{2m}|x_1| \right)} = f_{1}(x_1), \]
\[u_2 = \frac{2m}{2m-1} \log{ \left( \frac{2m-1}{2m}|x_2| \right)} = f_{2}(x_2), \]
\[u_3 = \frac{2m}{2m-1} \log{ \left( \frac{2m-1}{2m}|x_3| \right)} = f_{3}(x_3), \]
\[u_4 = -u_{1}-u_{2}-u_{3} = -\frac{2m}{2m-1} \log{ \left( \frac{2m-1}{2m}|x_4| \right)} = f_{4}(x_4). \]
So, as a separable hypersurface, we get
\[-\log{|x_1|}+\log{|x_2|}+\log{|x_3|}-\log{|x_4|} = 0, \]
which is equivalent to that
\[\frac{x_{2}x_{3}} {x_{1}x_{4}} = \pm 1. \]
We can show its minimality also directly.

\vspace{2mm}

Graduate School of Science and Technology 

Hirosaki University 

Hirosaki 036-8561, Japan 

E-mail: sakaki@hirosaki-u.ac.jp 

\end{document}